# A Solution for Large Scale Optimization Problems Based on Gravitational Search Algorithm

Somayeh Seifi Shalamzari[1], Mojtaba Banifakhr[2]
1- Shahid Bahonar University of Kerman: Kerman, Iran, somayeh.seifi1361@gmail.com.
2- Yazd University: Yazd, Iran, banifakhr.mojtaba@stu.yazd.ac.ir.

*Abstract*— One of the challenges in optimization of high dimensional problems is finding appropriate solutions in a way that are as close as possible to the global optima.. In this regard, one of the most common phenomena that occurs is the curse of dimensionality in which a large scale feature space generates more parameters that need to be estimated. Heuristic algorithms, such as Gravitational Search Algorithm (GSA), are among the tools proposed for optimizing large-scale problems, but in this case, they cannot solve the problem on their own. This paper proposes a novel method for optimizing large scale problems by improving the gravitational search algorithm's performance. In order to increase the efficiency of the gravitational search algorithm in solving large scale problems, the proposed method combines this algorithm with the cooperative-coevolution methods. For the evaluation of the performance of the proposed algorithm, we consider two approaches. In the first approach, the proposed algorithm is compared with the gravitational search algorithm, and in the second approach, it is compared with some of the most significant research in this field. In the first approach, our method was able to improve the performance of the original gravitational algorithm to solve large scale problems, and in the second one, the results indicate more favorable performance, in some benchmark functions, compared with other cooperative methods.

## 1. Introduction

Meta Heuristic algorithms are a branch of optimization algorithms that are inspired by nature and have the capability to search high dimensional spaces. During the last decades, different types of meta heuristic algorithms have been introduced. Some of these algorithms are: Artificial Bee Colony (ABC) [1], Ant Colony Optimization (ACO) [2], Particle Swarm Optimization (PSO) [3], an improved Shuffled Frog-leaping Algorithm (SFL) [4], Grenade Explosion Method (GME) [5], Artificial Immune System [6], Gravitational Search Algorithm (GSA) [7], and Whale Optimization Algorithm [8].

The gravitational search algorithm is one of the newest heuristic algorithms proposed by Rashedi and Nezam Abadipour in recent years, which is inspired by the concepts of mass and gravity by simulating Newton's laws of gravity and motion. This algorithm's features include simple implementation, the generality of the concept, and its ability to achieve global optimality, which have made it successful [6].

This algorithm is known as one of the most powerful and reliable search tools for the best solution, but as optimization problems become more complex and larger, researchers face new challenges in responding to the increasing computational needs of using this algorithm.

Optimization is one of the most important issues in the field of Computer Science, Artificial Intelligence, Soft Computing, and Machine Vision. In optimization, looking for the best solution in a reasonable time is desired [3-9]. Evolutionary algorithms have successfully solved most of the optimization problems, although their performance is weak in tackling high dimensional problems. Since the search space of a problem increases exponentially with the size of the problem, most of the evolutionary algorithms suffer from the "Curse of dimensionality." In this condition, the evolutionary algorithm is stuck in a local optimum and cannot find the global optimum [10].

In high dimensional problems, everything becomes much harder when we are facing with non-separable problems and interacting variables. For handling the obstacles related to the high dimensional problems, two categories of methods have been introduced. The first category is

Cooperative Co-evolutionary methods, and the second one is non Cooperative Coevolutionary methods. Cooperative Co-evolutionary (CC) approach is inspired by the divide and conquer method [10]. CC is a combination method that consists of two phases; the first phase is the grouping stage, and the second phase is the optimization stage. In the grouping stage, a high dimensional problem is decomposed into subcomponents and in the optimization stage, each subcomponent is optimized by an evolutionary algorithm. Since the most important stage in the CC method is the grouping stage, a suitable grouping strategy will guide the problem to the global optimum. The best grouping method is one which puts the interacting variables in the same subcomponents.

Formerly, CC has been applied to many evolutionary algorithms with different grouping methods to solve large scale problems. For the first time, CCGA, a combination of the CC method and the genetic algorithm, was introduced by Potter and Jong [10]. In their method, a problem with $n$ variables was decomposed to $n$ subgroups with one variable in each subgroup. F. Vandenberg applied CC to the Particle Swarm Optimization (PSO) algorithm [11]. In the grouping stage, the problems were decomposed to $k$ subcomponents with $s$ variables in each subcomponent. Liu applied CC to Fast Evolutionary Programming (FEP), but this algorithm was trapped in a local optimum in one of the non-separable functions [12]. Yang introduced a new grouping strategy, named Random Grouping, and applied it to the Differential Evolution (DE) [13]. Random Grouping increased the possibility of laying two interacting variables in the same group. However, when the number of interacting variables increased, the possibility decreased. As a result, Omidvar introduced Frequently Random Grouping, which increased this possibility [14]. The performance of Frequent Random Grouping deteriorated when the number of interacting variables went beyond 5 variables. To handle it, Omidvar presented Delta Grouping [15] but it showed weak performance where there were more than one non-separable subcomponents. To solve this issue, Omidvar introduced a new grouping strategy, named Differential Grouping, and applied CC with Differential Grouping to the Differential Evolution [16]. The differential grouping had a great performance to handle interacting variables.

In [17], EADE, which is a non-Cooperative Co-evolution method, has been introduced. This method is a modified version of the deferential evaluation algorithm (DE). However, in EADE, the mutation and crossover operators have been improved to be suitable for tackling large scale problems.

In this paper, we applied differential grouping, as a grouping strategy in the CC method, to the GSA to improve its performance in facing large scale problems. Since, thus far, the Gravitational Search Algorithm has not been applied to large scale optimization problems, the proposed method is innovative.

The remainder of this paper is arranged as follows. In section II the standard gravitational search algorithm is introduced. In section III the proposed method for tackling high dimensional problems is given. Section IV presents the results and comparison with the standard GSA. In the results, the comparison of our method with several Cooperative Co-evolutionary methods and one of the non Cooperative methods is presented. Finally, section V concludes the paper.

## 2. Research Foundations

As mentioned in the previous section, optimization is one of the most important issues for researchers in various fields. The structures proposed based on the evolutionary algorithms have good results in many optimization problems, but these results weaken as the number of dimensions increases. Today, with the advancement of various sciences and technologies, due to the larger scale of the problems and search space, the need for large scale optimization is arised.

It should be noted that when the search space grows, most optimization problems suffer from the curse of dimensionality, which deteriorates the performance of these methods [3-8]. One of the important factors that challenge large scale problems is that the increasing number of variables leads to an exponential increase in the problem space. In addition, characteristics may change as the scale of the problem space increases. Another challenge is the costly evaluation of large scale problems, and the last challenge is the existence of an interaction between variables.

## 2.1. Cooperative Co-evolutionary (CC) Method

The Cooperative Co-evolutionary method is a new method based on and inspired by the divide-and-conquer method, in which a large scale problem is divided into smaller scale subcomponents, and each subcomponent is optimized by a specific evolutionary algorithm (EA). In the CC method, the kind of grouping method we choose has a great impact on the optimization's performance. The better grouping we provide, the better the problem can be optimized and the probability of placing interacting variables in a group can be increased.

The hybrid optimization method with the CC method consists of two stages. In the first stage, grouping, the problem is grouped by applying one of the grouping methods, and thus subcomponents are formed. Then in the second stage, optimization, each of these subcomponents is optimized with a specific EA.

In the Gravitational Search Algorithm, optimization is performed with the help of the laws of gravity and motion in a discrete-time artificial system [6]. The system environment is the same as the scope of the problem definition. According to the law of gravitation, every mass perceives other masses' location and status through the law of gravitational attraction. Therefore, this force can be used as a tool to exchange information in space. The amount of masses is determined based on the objective function. The gravitational search algorithm is described in two general steps: in the first step, a discrete-time artificial system is formed in the problem environment, initial positioning of masses is determined, governing rules are defined, and parameters are set, and in the second step, the passage of time, the motion of masses, and the updating of parameters happen until it reaches the stop time.

## 2.2 Gravitational Search Algorithm

Gravitational Search Algorithm, which is a population-based optimization algorithm, is inspired by Newton's laws of gravity and motion [6]. Referred to [6], the mass of every agent is computed using the fitness of the population (for a minimization problem), as shown in Eq (1,2):

$$(1) \quad q_i(t) = \frac{fit_i(t) - worst(t)}{best(t) - worst(t)}$$

$$(2) \quad M_i(t) = \frac{q_i(t)}{\sum_{j=1}^{N} q_j(t)}$$

Where $N$, $M_i(t)$ and $fit_i(t)$ represent the size of the population, the mass, and fitness of the $i^{th}$ agent, resepctively. Worst (t) and best (t) are calculated (for a minimization problem) based on Eq (3,4).

$$(3) \quad best(t) = \min_{j \in \{1,\ldots,N\}} fit_j(t)$$

$$(4) \quad worst(t) = \min_{j \in \{1,\ldots,N\}} fit_j(t)$$

Total forces of heavier masses applied to an agent based on the laws of gravity and motion [18] in Eq (5-7) is used to compute the acceleration of an agent as shown in Eq (6). The next velocity, which is a fraction of the current velocity added to the acceleration is calculated in Eq (7). The next position of an agent will be calculated based on Eq (8):

$$(5) \quad F_{ij}^d(t) = G(t) \frac{M_{pi}(t) M_{aj}(t)}{R_{ij}(t) + \varepsilon} \left( x_j^d(t) - x_i^d(t) \right)$$

$$(6) \quad R_{ij}(t) = \|X_i(t), X_j(t)\|_2$$

$$(7) \quad F_i^d(t) = \sum_{j \in kbest, j \neq i} rand_j F_{ij}^d(t)$$

$$(8) \quad a_i^d = \frac{F_i^d(t)}{M_{ii}(t)}$$

$$(9) \quad v_i^d(t+1) = rand_i \cdot v_i^d(t) + a_i^d(t)$$

$$(10) \quad x_i^d(t+1) = x_i^d(t) + v_i^d(t+1)$$

In these formulas, $v_i^d$, $a_i^d$ and $x_i^d$ respectively depict the acceleration, velocity, and position of the $i^{th}$ agent in $d^{th}$ dimension. $rand_i$ and $rand_j$ are random numbers in the interval [0,1]. Epsilon is a small number, n is the dimension of the search space and R is the Euclidean distance between two agents. Agents with the best fitness and the greatest masses are presented by $K_{best}$. $K_{best}$ is initialized by $K_0$ equal to N and decreases linearly to 1 by time.

G is initialized by $G_0$ and decreases exponentially to zero by time. The block diagram of the standard GSA is presented in figure 1.

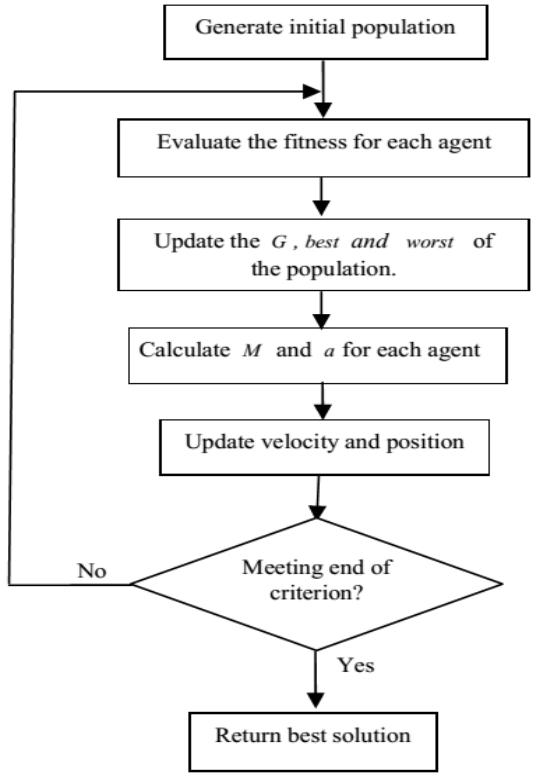

**Fig. 1.** Block diagram of the standard GSA

So far, several modified versions of GSA have been introduced such as Binary Gravitational Search Algorithm [19], Disruption: a new operator in gravitational search algorithm [20], a quantum version of gravitational search algorithm [21], and Black hole: a new operator for gravitational search algorithm [22].

### 3. Proposed Method

As mentioned earlier, the GSA algorithm was only suitable for small scale functions. Therefore, for larger-scale problems, an algorithm is needed that could cross local minima in addition to having a good convergence speed.

Like other cooperative co-evolutionary methods, the proposed method consists of two stages: in the first stage, the grouping is performed, and then in the next stage, optimization is conducted.

In the grouping stage, any of the grouping methods described can be used. As stated earlier, the differential grouping algorithm is the best grouping method, as in this algorithm, the interacting variables are most likely to fall in the same subgroup. Therefore, the proposed method uses the differential grouping algorithm in the grouping stage. In this method, first, the differential grouping algorithm is used to separate dimensions, and then the GSA algorithm, which provides an acceptable solution compared to other metahuristic algorithms, is used as the optimizer.

Accordingly, the fitness value calculated in the current algorithm has a significant advantage over the GSA. The difference between the proposed method and other cooperative co-evolutionary methods is that in this implemented algorithm, first, the optimizer function is different, and thus, the inputs and outputs are slightly different, and second, each time the optimization per column (dimension) is performed, the best value obtained is stored in the corresponding column and in a variable. According to the approach of using differential grouping, the calculated groups are independent of each other. Therefore, in the calculation, each column can be optimized independently. A noteworthy point in the proposed algorithm is that each time the algorithm is executed, the best member of the population is identified, and then a population is constructed in the optimizer function that has the size of the main population, and also all the rows of this population are made up of this best member. In the optimization phase, the subgroup to which the optimizer function is to be applied takes its place in this new population, and the optimizer function applies only to this subgroup.

the advantage of this proximity to the best member is that with consecutive iterations of the algorithm, all the population members move toward this best member. Using the above algorithm, although the computation time increases slightly, by looking at the convergence diagrams, one can see that this trend gets better in most diagrams. It should be noted that in some diagrams, convergence sometimes declines, but then it is observed that the previous decline is offset and moves toward improvement. As can be deduced from the pseudo-code presented in Fig. 2, the stated operation is reiterated by the number of cycles.

Algorithm 1

Begin

Initialization (N,dim,LB,UB,FEs)

    Groups←Differential_Grouping(Func_num,LB,UB,dim)

    pop←rand(N,dim)

    for i=1 to cycles

        best←min(func(pop))

        for j=1 to size(Groups)

            index←Groups[j]

            subpop←pop[:,index]

        [subpop,BestValue]←GSA(subpop,best,Fes)

            pop[:,index]←subpop

            Result[index]← BestValue

        endfor

    endfor

End

**Fig. 2.** Pseudo-code of the proposed algorithm (CCGSA-DG)

## 4. Simulations

In the previous section, a hybrid CCGSA-DG algorithm was introduced in which differential grouping is used as a proper grouping algorithm. In this section, this algorithm is compared with other algorithms. This comparison is performed in two stages. In the first stage, the algorithm is compared with the standard GSA. The benchmark functions available in [1] are used for this purpose. In the second stage, the proposed method is compared with the following algorithms: DECC-DG [16], MLCC [23], DECC-D [15], DECC-DML [15], and EADE [17].

### 4-1 Comparison of CCGSA-DG with the standard GSA

In this section, the benchmark functions available in [1] are employed for evaluation. There are 24 benchmark functions, among which 13 functions, F1 to F13, can have different number of variables, but other functions, F14 to F23, have a fixed number of variables. Functions $F_1$ to $F_7$ are unimodal, in which the convergence rate of the searching algorithm is more important than the final answer, and functions $F_8$ to $F_{13}$ are multimodal, in which number of minimums for the above functions increases exponentially as the number of dimensions increases, making the problem more complex. The list of unimodal and multimodal functions used in this experiment is given in tables 1 and 2. Since we tend to examine the existing algorithm's performance on large scale functions and compare it with other cooperative algorithms of this field, this algorithm is first applied only to some specific benchmark functions from [1], and then it is applied to all 2010 benchmark functions [3] for comparison with other cooperative algorithms. In this table, n, S, and $f_{opt}$ indicate the search space dimension, search space region, and the optimal value, respectively.

**Table 1.** Benchmark functions used to compare the suggested algorithm with GSA (unimodal functions) [1].

| Test Function | S | $f_{opt}$ |
|---|---|---|
| $F_1(X) = \sum_{i=1}^{n} x_i^2$ | $[-100,100]^n$ | 0 |
| $F_4(X) = \max_i\{|x_i|, 1 \leq i \leq n\}$ | $[-100,100]^n$ | 0 |
| $F_6(X) = \sum_{i=1}^{n}([x_i + .5])^2$ | $[-100,100]^n$ | 0 |

**Table 2.** Benchmark functions utilized to compare the proposed algorithm with GSA (multimodal functions) [1].

| Test Function | S | $f_{opt}$ |
|---|---|---|
| $F_8(X) = \sum_{i=1}^{n} -x_i \sin(\sqrt{|x_i|})$ | $[-500,500]^n$ | $-418.9829 \times n$ |
| $F_9(X) = \sum_{i=1}^{n}[x_i^2 - 10\cos(2\pi x_i) + 10]$ | $[-5.12,5.12]^n$ | 0 |
| $F_{10}(X) = -20\exp\left[-.2\sqrt{\frac{1}{n}\sum_{i=1}^{n}x_i^2}\right] - \exp\left[\frac{1}{n}\sum_{i=1}^{n}\cos(2\pi x_i)\right] + 20 + e$ | $[-32,32]^n$ | 0 |
| $F_{11}(X) = \frac{1}{4000}\sum_{i=1}^{n}x_i^2 - \prod_{i=1}^{n}\cos\left(\frac{x_i}{\sqrt{i}}\right) + 1$ | $[-600,600]^n$ | 0 |

Table 3. Performance of the proposed method compared to GSA according benchmark functions

| function | dimention | CCGSA-DG | GSA |
|---|---|---|---|
| $F_1$ | 30 | $5.62 \times 10^{-22}$ | $2.6 \times 10^{-17}$ |
| $F_1$ | 500 | $1.23 \times 10^{-8}$ | $3.2 \times 10^4$ |
| $F_1$ | 1000 | $5.5 \times 10^{-8}$ | $3.8 \times 10^4$ |
| $F_4$ | 30 | 3.12 | 1.6 |
| $F_4$ | 500 | 98.84 | 22.36 |
| $F_4$ | 1000 | 99.47 | 25.5 |
| $F_6$ | 30 | 0 | 0 |
| $F_6$ | 500 | 0 | 36561 |
| $F_6$ | 1000 | 0 | 81728 |
| $F_8$ | 30 | -12478.932 | $-2.5 \times 10^3$ |
| $F_8$ | 500 | -206525.81 | $-1.01 \times 10^4$ |
| $F_8$ | 1000 | -413157.27 | $-1.4 \times 10^4$ |
| $F_9$ | 30 | 0 | 21.8 |
| $F_9$ | 500 | 0.0068629 | $2.084 \times 10^3$ |
| $F_9$ | 1000 | 0.017156 | $5.6 \times 10^3$ |
| $F_{10}$ | 30 | 17.88 | $4.7 \times 10^{-9}$ |
| $F_{10}$ | 500 | 18.54 | 8.34 |
| $F_{10}$ | 1000 | 19.64 | 9.06 |
| $F_{11}$ | 30 | 0.49 | 16.6 |
| $F_{11}$ | 500 | 0.048484 | $8.37 \times 10^3$ |
| $F_{11}$ | 1000 | 0.068055 | $1.99 \times 10^4$ |

In the GSA algorithm applied in CCGSA-DG, the population size (N) and the maximum iterations were considered 50 and 500, respectively. In equation G used in GSA, the constant $G_0$ is considered 100 and $\propto$ is considered equal to 20.

$$(11) \quad G = G_o \exp(-\frac{\propto . t}{t_{max}})$$

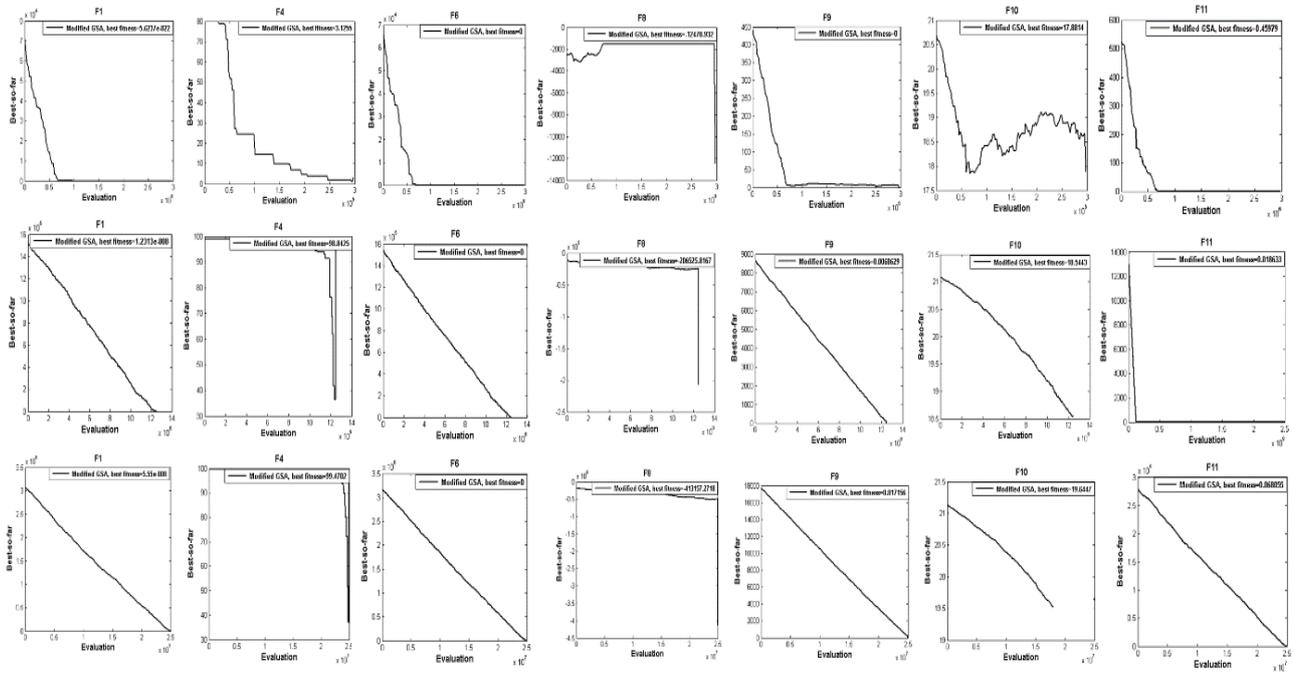

**Figure 3.** Convergence diagram for desire functions in three dimension 30(up) 50(midle) 1000(down)

In the simulation conducted in this section, the maximum of the cycles that optimization is applied to subgroups is considered equal to 20.

Figure 3 shows the results obtained from applying the proposed CCGSA-DG algorithm on functions with 30, 500, and 1000 dimensions. As the obtained convergence process of figures and optimal value show, the resulting algorithm has provided favorable results. As can be concluded from the above figure, convergence graphs for 500 dimensions, which is a lot more than that of the original GSA, have presented acceptable results and convergence.

As we have stated, the gravitational algorithm has not been applied to functions with high dimensions. According to the graphs, the suggested algorithm has represented good results on functions with 1000 dimensions.

## 4-2 Comparing the CCGSA-DG algorithm with existing algorithms in the field of CC

In this section, the suggested algorithm is compared with algorithms DECC-DG, MLCC, DECC-DML, DECC-D, and EADE. The benchmark functions employed in this section are 2010 benchmark functions [10] written for functions with high dimensions. The proposed algorithm was applied to these functions with 1000 dimensions. The 2010 benchmark functions are usually used for problems with high dimensions. These functions are designed in such a way that they solve four types of large scale problems. In the first group, separable functions are used as benchmark functions, which include functions $F_1$ to $F_3$. In the second group, the partial separable functions are used as benchmark functions. In these functions, there is a non-separable subcomponent with 50 variables and a separable group consisting of 950 variables. Functions $F_4$ to $F_8$ belong to this group. In the third group, some functions include 10 non-separable subcomponents, each containing 50 variables, and a separable subcomponent with 500 variables. Functions $F_9$ to $F_{13}$ belong to this group. The fourth group contains functions that include 20 non-separable subcomponents, each with 50 variables. Functions $F_{14}$ to $F_{18}$ are in this group. The last group includes functions that are completely non-separable and contain functions $F_{19}$ to $F_{20}$.

In all the functions, the search space dimensions are considered to be equal to 1000, and $p$ is a random permutation of set $\{1, 2, ..., D\}$. In Table 4, the performance of the proposed method compared to other methods considered in the second approach can be observed by taking into account the above-mentioned benchmark function.

**Table 4.** Performance of the proposed method compared to other methods according benchmark functions

| function | CCGSA-DG | DECC-DG | MLCC | DECC-D | DECC-DML | EADE |
|---|---|---|---|---|---|---|
| $F_1$ | $1.1 \times 10^8$ | $5.47 \times 10^3$ | $1.53 \times 10^{-27}$ | $1.01 \times 10^{-24}$ | $1.93 \times 10^{-25}$ | $4.7 \times 10^{-22}$ |
| $F_2$ | $6.04$ | $4.39 \times 10^3$ | $5.57 \times 10^{-1}$ | $2.99 \times 10^2$ | $2.17 \times 10^2$ | $4.16 \times 10^2$ |
| $F_3$ | $1.26$ | $1.67 \times 10^1$ | $9.88 \times 10^{-12}$ | $1.81 \times 10^{-13}$ | $1.18 \times 10^{-13}$ | $6.25 \times 10^{-14}$ |
| $F_4$ | $8.2 \times 10^{14}$ | $4.7 \times 10^{12}$ | $9.6 \times 10^{12}$ | $3.84 \times 10^8$ | $3.5 \times 10^{12}$ | $1.08 \times 10^{11}$ |
| $F_5$ | $3.5 \times 10^8$ | $1.55 \times 10^8$ | $3.84 \times 10^8$ | $4.16 \times 10^8$ | $2.98 \times 10^8$ | $8.79 \times 10^7$ |
| $F_6$ | $1.9 \times 10^7$ | $1.64 \times 10^1$ | $1.62 \times 10^7$ | $1.36 \times 10^7$ | $7.93 \times 10^5$ | $1.9 \times 10^1$ |
| $F_7$ | $2.8 \times 10^{12}$ | $1.16 \times 10^4$ | $6.89 \times 10^5$ | $6.58 \times 10^7$ | $1.39 \times 10^8$ | $2.11 \times 10^1$ |
| $F_8$ | $2.17 \times 10^{16}$ | $3.04 \times 10^7$ | $4.38 \times 10^7$ | $5.39 \times 10^7$ | $3.46 \times 10^7$ | $2.26 \times 10^{-4}$ |
| $F_9$ | $6.5 \times 10^{10}$ | $5.96 \times 10^7$ | $1.23 \times 10^8$ | $6.19 \times 10^7$ | $5.92 \times 10^7$ | $3.67 \times 10^7$ |
| $F_{10}$ | $6.1 \times 10^3$ | $4.52 \times 10^3$ | $3.43 \times 10^3$ | $1.16 \times 10^4$ | $1.25 \times 10^4$ | $2.62 \times 10^3$ |
| $F_{11}$ | $2.06 \times 10^2$ | $1.03 \times 10$ | $1.98 \times 10^2$ | $4.76 \times 10^1$ | $1.8 \times 10^{-13}$ | $1.14 \times 10^2$ |
| $F_{12}$ | $6.1 \times 10^6$ | $2.52 \times 10^3$ | $3.49 \times 10^4$ | $1.53 \times 10^5$ | $3.79 \times 10^6$ | $2.8 \times 10^4$ |
| $F_{13}$ | $2.7 \times 10^{11}$ | $4.54 \times 10^6$ | $2.08 \times 10^3$ | $9.87 \times 10^2$ | $1.14 \times 10^3$ | $1.01 \times 10^3$ |
| $F_{14}$ | $2.5 \times 10^{10}$ | $3.41 \times 10^8$ | $3.16 \times 10^8$ | $1.98 \times 10^8$ | $3.38 \times 10^8$ | $1.46 \times 10^8$ |
| $F_{15}$ | $4.2 \times 10^3$ | $5.88 \times 10^3$ | $7.11 \times 10^3$ | $1.53 \times 10^4$ | $5.87 \times 10^3$ | $3.18 \times 10^3$ |
| $F_{16}$ | $2.94 \times 10^2$ | $7.39 \times 10^{-13}$ | $3.76 \times 10^2$ | $1.88 \times 10^2$ | $2.47 \times 10^{-13}$ | $3 \times 10^2$ |
| $F_{17}$ | $1.1 \times 10^6$ | $4.01 \times 10^4$ | $1.59 \times 10^5$ | $9.03 \times 10^5$ | $3.19 \times 10^4$ | $1.52 \times 10^5$ |
| $F_{18}$ | $8.5 \times 10^{11}$ | $1.11 \times 10^{10}$ | $7.09 \times 10^3$ | $2.12 \times 10^3$ | $1.17 \times 10^3$ | $2.26 \times 10^3$ |
| $F_{19}$ | $3 \times 10^6$ | $1.47 \times 10^6$ | $1.36 \times 10^6$ | $1.33 \times 10^7$ | $1.74 \times 10^6$ | $1.29 \times 10^6$ |
| $F_{20}$ | $8.9 \times 10^{11}$ | $4.47 \times 10^7$ | $2.05 \times 10^3$ | $2.05 \times 10^3$ | $4.14 \times 10^3$ | $2.1 \times 10^3$ |

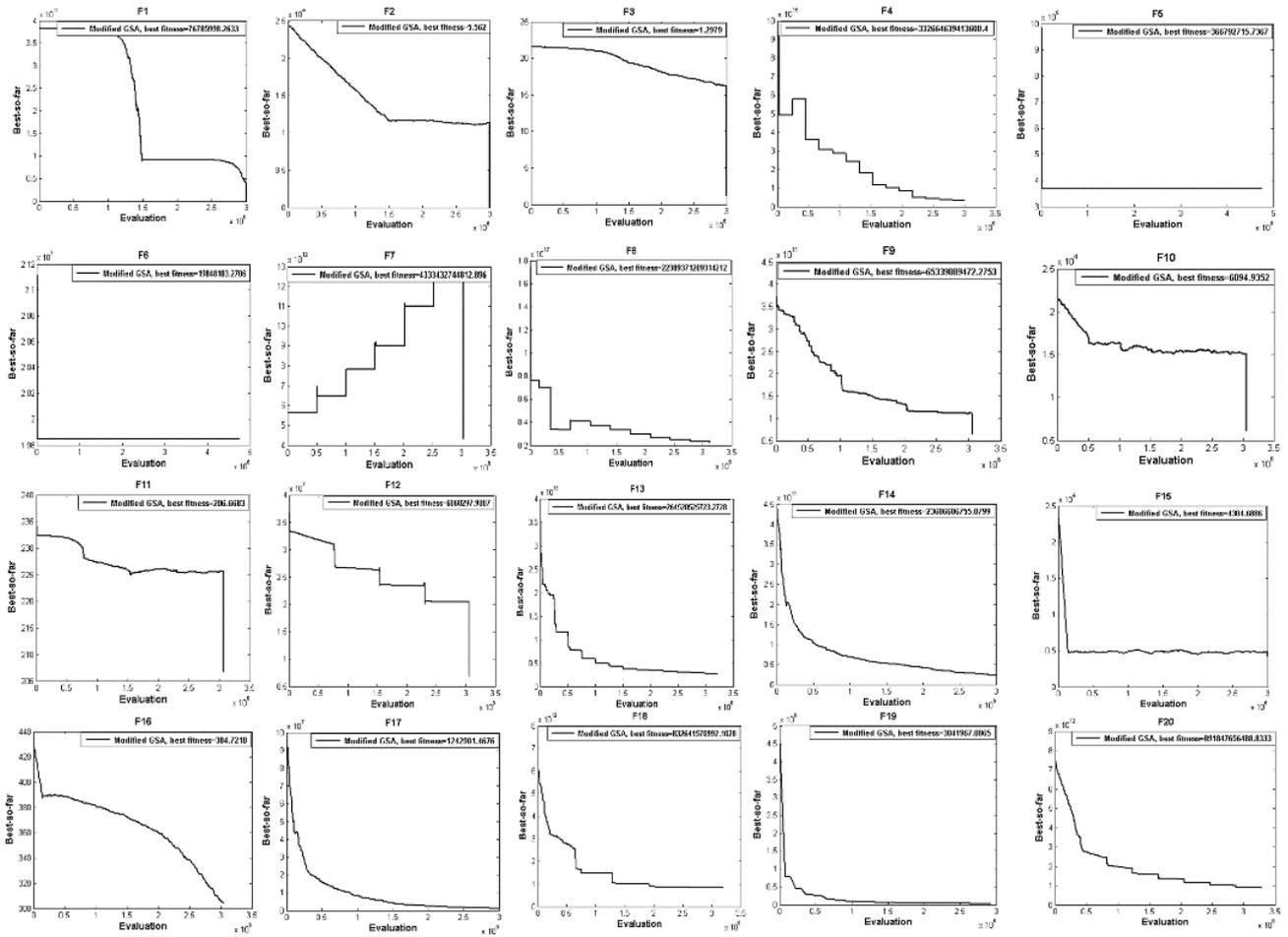

**Figure 4**. Convergence graphs for the second approach

The convergence graphs associated with existing functions are given in figure 4. In the top, middle, and bottom figures, the convergence graphs are shown for dimension equal to 1000.

## 4-3 Results

For the first considered approach in this paper, the *GSA's* ability to solve problems with high dimensions is improved by applying *CC* to the *GSA* algorithm and using differential grouping as a proper grouping strategy. In this approach, the maximum number of cycles in which optimization is applied to subgroups is considered to be 20 and tested on functions with 30, 500, and 1000 dimensions. Since the original *GSA* algorithm was only investigated when applied to functions with up to 30 dimensions, a comparison of the suggested *GSA* algorithm with the original *GSA* shows that this algorithm achieved promising results in functions with 30 dimensions compared with the original *GSA*. This algorithm was tested on functions with 500 and 1000 dimensions as well, and the results have shown that it performed well on these dimensions. As can be deduced from the convergence process of figures and the obtained optimal value, the resulting algorithm has achieved good results in most cases.

In the second approach, according to the results, this algorithm has a good performance in nonseparable functions except for $F_1$. These functions are functions of the first group, including $F_1$, $F_2$, and $F_3$. As can be concluded from the table's values, the proposed algorithm was able to achieve an almost optimal point. The maximum number of cycles, the population size (N), and the maximum number of algorithm evaluations are considered 20, 50, and 3000000, respectively. The values shown in Table 4 are obtained from 25 separate executions. As Table 4 represents, in some functions, the current algorithm had better performance compared to some of the existing algorithms, and in some other, it did not perform favoravly. The notable point is that if we consider convergence graphs, benchmark functions have shown a good convergence in most cases. As can be deduced from graphs, in most cases, the graph is decremental, and if we could increase the number of iterations, we could accomplish more optimal values. However, considering that these functions are large scale, this is time-consuming and requires a high-speed, powerful computation system. It is highly probable that such a system improves the answers.

In addition, the suggested algorithm had a better performance in function $F_2$ compared with algorithms DECC-DG, DECC-D, and DECC-DML. Since $F_2$ is a separable function, the proposed method had quite good results for the separable functions. In function $F_3$, which is also a separable function, our method performed better than algorithms MLCC and DECC-D. Moreover, in function $F_{10}$, which is a partial seperable function, it had a better performance than DECC-D and DECC-DML. The proposed method achieved better results than all other available algorithms except EADE in function $F_{15}$, which is a function with 20 non-separable subcomponents. The proposed algorithm performed better than MLCC , DECC-D and EADE in function $F_{16}$, which is a function with 20 non-separable subcomponents and it has shown better results than DECC-D in $F_{19}$.

## 5- Conclusion

The Swarm intelligence algorithms were successful in solving most of the optimization problems. However, these algorithms faced the curse of dimensionality when confronting high dimension problems, which weakened their performance. One of the most popular methods proposed to solve large scale problems is using the cooperative-coevolutionary (CC) methods inspired by the divide and conquer approach. In this approach, the problem space is broken into subcomponents, each of which is optimized separately using an EA. Thus, instead of dealing with a vast space, we have subspaces with smaller dimensions that are much easier to optimize as compared with the original problem. The most important factor in the CC method is choosing the proper decomposition method with the highest accuracy percentage in placing the interacting variables in subgroups. One of the best grouping methods with the highest accuracy is the differential grouping method. In the suggested method in this paper, a problem with high dimensions is broken down into subcomponents, and a GSA optimizer algorithm is applied to each one. The proposed method was examined with two approaches. In the first approach, there was a tangible improvement in results for all benchmark functions, and in the second one, the performance of some functions has improved.

## REFERENCES

[1] D. Karaboga, B. Basturk, J.Glob, "A Powerful and Efficient Algorithm for Numerical Function Optimization: Artificial Bee Colony (ABC) Algorithm Optim", pp.459-471, (2007).

[2] M. Dorigo and T. Stutzle, "Ant Colony Optimization", MIT Press, ISBN 0-262-04219-3, 2004.

[3] J. Kennedy, R.C Eberhart, "Particle Swarm Optimization". In: Proceedings of the 1995 IEEE, International Conference on Neural Networks, IEEE Service Center, Piscat away, PP. 1942-1948, 1999.

[4] X. Li, J. Luo, M-R. Chen, N. Wang, "An Improved Shuffled Frog-leaping Algorithm with Extremal Optimization for Continuous Optimization", Information Sciences, Vol.192, pp. 143-151, 2010.

[5] A. Ahrari, A. A. Atia, "Grenade Explosion Method-A Novel Tool for Optimization of Multimodal Functions", Applied Soft Computing, pp.1132-1140, 2010.

[6] D. Castro, L.N. Von, F.J. Zuben, "Artificial Immune Systems". Part I. Basic Theory And Applications, 1999.

[7] Esmat Rashedi, Hossein Nezamabadi-pour, Saeid Saryazdi, "A Gravitational Search Algorithm", Information sciences, pp.2232-2248, 2009.

[8] Seyedali Mirjalali, Andrew Lewis, " The Whale Optimization Algorithm", Advances in Engineering Software, pp.51-67, 2016.

[9] Hossein Nezamabadipour, "Inheritance Algorithm, Basic and Advanced Concepts", First Edition, 2010.

[10] Ke Tang, Xiaodong Li, P.N.Suganthan, Zhenyn Yang, Thomas Weise, "Benchmark Functions for the CEC2010 Special Session and Competition on Large –Scale Global Optimization", pp. 1-23, 2009.

[11] F. van den bergh and A. Engelbrecht, "A Cooperative Approach to Particle Swarm Optimization", IEEE Trans. Evol. Comut, Vol.8, No. 3, pp. 225-239, 2004.

[12] Y. Liu, X. Yao, Q. Zhao, and T. Higuchi, "Scaling up Fast Evolutionary Programming with Cooperative Coevolution", In proceedings of 2001 Congress on Evolutionary Computation, pp. 1101-1108, 2001.

[13] Z. Yang, k. Tang, and X. Yao, "Large Scale Evolutionary Optimization Using Cooperative Coevolution", Information Sciences, Vol.178, No.15, pp.2985-2999, 2008.

[14] M. N. Omidvar and X. Li, "Cooperative Co-evolution for Large Scale Optimization Through More Frequent Random Grouping", IEEE, pp.1-8, 2010.

[15] M. N. Omidvar and X. Li, "Cooperative Co-evolution with Delta Grouping for Large Scale Non-separable Function Optimization", IEEE, pp.1-8, 2010.

[16] M. N. Omidvar, X. Li, Y. Mei and X. Yao, "Cooperative Co-evolution with Differential Grouping for Large Scale Optimization", IEEE transection on evolutionary computation, Vol. X, No.X, pp. 1-16, 2013.

[17] Ali Wagdy Mohamed, "Solving large scale global optimization problems using enhanced adaptive differential evolution algorithm", 2017

[18] D. Holiday, R. resnick and J. Walker, "Fundamentals of Physics", John Wiley and sons, ISNB 0470469080, 1993.

[19] Esmat Rashedi, Hossein Nezam Abadi Pour, Saeid Saryazdi, "BGSA:binary gravitational search algorithm",Natural Computing, Vol.9, pp.727-745, 2010

[20] S Sarfarazi, H Nezamabadi Pour, S Saryazdi, "Disryption: a new operator in gravitational search algorithm", Scientia Iranica, Vol.18, 2011

[21] Mohadeseh Soleimanpour-Moghadam, Hossein Nezamabadi-Pour, Malihe M Farsangi,"A quantum inspired gravitational search algorithm for numerical function optimization", Information Sciences, Vol.267, pp.83-100, 2014.

[22] Mohammad Doraghinejad, Hossein Nezamabadi-Pour, "Black hole: a new operator for gravitational search algorithm", International journal of Computational Intelligence System, Vol.7, pp.809-826, 2014

[23] Z. Yang, K. Tang, and X. Yao, "Multilevel Cooperative Coevolution for Large Scale Optimization", in proc. IEEE Congr. Evol. Comput, jun. pp. 1663-1670, 2008.